\documentstyle[psfig]{siamltex}

\title{
ASYMPTOTIC-NUMERICAL STUDY OF SUPERSENSITIVITY
FOR GENERALIZED BURGERS EQUATIONS}

\author{
   Marc Garbey\thanks{
      CDCSP-ISTIL,
      Universit\'e Claude Bernard Lyon~1,
      69622~Villeurbanne cedex, France
      (garbey@cdcsp.univ-lyon1.fr).
      Supported by the
      Fondation Cromey-le-Bas
      under contract BS-95-2.
   }
   \and Hans G.~Kaper\thanks{
      Mathematics and Computer Science Division,
      Argonne National Laboratory,
      Argonne, IL 60439-4801, USA
      (kaper@mcs.anl.gov).
      Supported by the
      Mathematical, Information, and Computational Sciences Division
      subprogram of Advanced Scientific Computing Research,
      U.S.\ Department of Energy,
      under contract W-31-109-Eng-38.
   }
}

\begin{document}

\maketitle

\begin{abstract}
This article addresses some asymptotic and
numerical issues related to the solution of
Burgers' equation,
$-\varepsilon u_{xx} + u_t + u u_x = 0$ on $(-1,1)$,
subject to the boundary conditions
$u(-1) = 1 + \delta$, $u(1) = -1$,
and its generalization to two dimensions,
$-\varepsilon \Delta u + u_t + u u_x + u u_y = 0$
on $(-1,1) \times (-\pi, \pi)$,
subject to the boundary conditions
$u|_{x=1} = 1 + \delta$, $u|_{x=-1} = -1$,
with $2\pi$ periodicity in $y$.
The perturbation parameters $\delta$ and $\varepsilon$
are arbitrarily small positive and independent;
when they approach 0, they satisfy
the asymptotic order relation
$\delta = O_s ({\rm e}^{-a/\varepsilon})$
for some constant $a \in (0,1)$.

The solutions of these convection-dominated
viscous conservation laws exhibit
a transition layer in the interior
of the domain,
whose position as $t\to\infty$
is supersensitive to the boundary perturbation.
Algorithms are presented for the computation
of the position of the transition layer
at steady state.
The algorithms generalize to viscous conservation
laws with a convex nonlinearity and are scalable
in a parallel computing environment.
\end{abstract}

\pagestyle{myheadings}
\thispagestyle{plain}
\markboth{MARC GARBEY AND HANS G.~KAPER}
         {SUPERSENSITIVITY FOR GENERALIZED BURGERS EQUATIONS}



\section{Introduction}
In this article we address some asymptotic and
numerical issues related to the solution of
Burgers' equation,
\begin{equation}
   - \varepsilon u_{xx}
   + u_t
   + u u_x
   = 0
   \quad\mbox{on } (-1, 1) ,
   \quad
   u(-1) = 1 + \delta ,
   \quad
   u(1) = - 1 ,
   \label{B1}
\end{equation}
and its generalization to two dimensions,
\begin{equation}
   - \varepsilon \Delta u
   + u_t
   + u u_x
   + \beta u u_y
   = 0
   \quad\mbox{on } (-1, 1) \times (-\pi, \pi) ,
   \quad
   u|_{x=-1} = 1 + \delta ,
   \quad
   u|_{x=1} = - 1 .
   \label{B2}
\end{equation}
In the latter case, we assume periodicity
(period $2\pi$) in the second coordinate ($y$).
The perturbation parameters $\delta$ and $\varepsilon$
are arbitrarily small positive;
they are independent, but when they approach 0,
they satisfy the asymptotic order relation
\begin{equation}
   \delta = O_s( {\rm e}^{-a/\varepsilon} )
   \quad
   \mbox{as } \delta, \varepsilon \downarrow 0 ,
   \label{delta-eps}
\end{equation}
for some constant $a \in (0,1)$
which does not depend on $\delta$ or $\varepsilon$.
This asymptotic relation implies that
$\delta$ is transcendentally small
(in the sense of asymptotic analysis)
compared with $\varepsilon$, but
$\delta$ dominates ${\rm e}^{-1/\varepsilon}$
as $\varepsilon \downarrow 0$.
(See~\cite{e79,s85,om91,kc96,djf96}
for definitions and basic concepts of
asymptotic analysis.)

If $\varepsilon = 0$, the solution of~(\ref{B1})
develops a shock (discontinuity)
in finite time, even when the initial data
are smooth~\cite{l73,t92}.
The perturbation introduced by a
nonzero $\varepsilon$ models the presence
of viscosity, which tends to smooth
the discontinuity~\cite{h50,f88}.
Instead of a shock, one has a
transition layer---a region of
rapid variation, which extends
over a distance $O(\varepsilon)$ as
$\varepsilon \downarrow 0$.
The position of the transition layer varies
with time, and its eventual location at
steady state is extremely sensitive
to the boundary data.
In fact, even the transcendentally small
perturbation $\delta$ leads to
a measurable (that is, order one)
effect on the eventual location
of the transition layer.
This phenomenon,
known as \textit{supersensitivity},
was first observed by Lorentz~\cite{l81}.
It has been studied extensively for Burgers' equation
and more general viscous conservation laws
in one dimension by
Kreiss and Kreiss~\cite{kk86},
Kreiss~\cite{k91},
Laforgue and O'Malley~\cite{lom93a,lom93b,lom94,lom95a,lom95b,lom99},
and Reyna and Ward~\cite{rw95a,rw95b,wr95}.

An example from combustion theory shows that
supersensitivity is of more than mathematical significance.
A simple model of flame propagation in gaseous fuels
involves a system of two coupled convection-diffusion
equations, one for the temperature of the mixture,
another for the concentration of the reaction-limiting
component in the mixture~\cite[\S~3.2]{bl83}.
If one ignores exponentially small perturbations
in the data, one finds that the Lewis number ${\cal L}$,
which is a measure for the ratio of
heat and mass transfer in the mixture, 
has no effect on the location of the combustion front.
Yet, numerical computations show that this location
is very sensitive---in fact, supersensitive---to
the value of ${\cal L}$.

Although the phenomenon of supersensitivity
is fairly well understood theoretically,
at least for one-dimensional problems,
the numerical solution of such problems
still poses formidable challenges,
especially in more than one dimension.
The methods that have been proposed
in the numerical literature for
singularly perturbed boundary-value
problems (see, for example, \cite{rst96})
tend to focus on uniform approximations
or finite elements with special features, 
not on the supersensitive dependence of
the transition layer on the boundary data.
On the other hand, the algorithms
we propose are designed specifically
to capture the phenomenon of supersensitivity.
They use the fact that the solution
approaches a certain profile
as the perturbation parameters
approach zero and focus on the
computation of the corrections.

Our ultimate goal is to develop algorithms
for multidimensional problems that are,
first of all, suitable for long-time
integration,
so stable steady states can be computed
with confidence;
second, extremely accurate in space,
so the eventual location of transition layers
can be predicted with accuracy; and third,
scalable in a multiprocessing environment,
so large-scale problems can be solved in a
reasonable length of time.
Although we discuss only Burgers' equation
and its generalization to two dimensions,
the algorithms are not restricted by
the special form of the nonlinearity.

In \S~2 we consider Burgers' equation.
We propose a simple algorithm that effectively
captures the supersensitive location of the
transition layer at steady state.
We stress the importance of the regions
outside the transition layer,
where the solution does not yet deviate
appreciably from the boundary values.
In \S~3 we address the generalized
Burgers equation in two dimensions.
We show through a formal asymptotic analysis
that the location of the transition layer
may vary in the direction of periodicity ($y$),
but the transition layer is essentially flat,
and only its average position (averaged over $y$)
depends supersensitively on the small parameters.
We then develop an algorithm that effectively
approximates the transition layer.

\section{One-Dimensional Case}
We begin by considering Eq.~(\ref{B1}),
\begin{equation}
   - \varepsilon u_{xx}
   + u_t
   + u u_x
   = 0
   \quad\mbox{on } (-1, 1) ,
   \quad
   u(-1) = 1 + \delta ,
   \quad
   u(1) = - 1 .
   \label{B}
\end{equation}
As shown by
Laforgue and O'Malley~\cite{lom95a},
the solution approaches a certain profile
function as $\varepsilon \downarrow 0$,
\begin{equation}
  u(x, t)
  =
  \tanh \eta
  +
  {\rm e}^{-a/\varepsilon} u_1 (\eta, \sigma)
  + \ldots \,,
  \label{tanh}
\end{equation}
where
\begin{equation}
  \eta = \frac{x - x^*(\sigma)}{\varepsilon} ,
  \quad
  \sigma = t {\rm e}^{-a/\varepsilon} .
\end{equation}
The hyperbolic tangent incorporates
a transition layer centered at $x^*$,
which connects the limiting values
$\pm 1$ at $\mp\infty$.
Note that these limiting values are
transcendentally close to the prescribed
boundary values $1+\delta$ and $-1$
of $u$ at $-1$ and $1$.
The position of the transition layer
varies on a transcendentally slow
time scale;
if $\delta = 2b {\rm e}^{-a/\varepsilon}$,
its limit as $\sigma \to \infty$ is
\begin{equation}
  x^*_{\rm as} = 1 - a + \varepsilon \ln b .
  \label{de-as}
\end{equation}
The important points to observe are that, first,
the solution $u$ approaches a certain
profile function as $\varepsilon \downarrow 0$;
second, the accurate determination of
the steady-state position of
the transition layer requires
long-term integration;
and, third, a transcendentally small
perturbation of the boundary data
has a measurable effect on the
location of the transition layer.

The asymptotic analysis has been generalized
to more general nonlinearities
by Laforgue and O'Malley~\cite{lom95b,lom99}
and Reyna and Ward~\cite{rw95a},
with similar  results.
One finds a limiting profile,
which generalizes the hyperbolic tangent function,
and a transition layer which moves
on a transcendentally slow time scale
to a steady-state position.
This position depends supersensitively
on the boundary perturbation.
Whenever this situation arises,
appropriate variants of the following
algorithms can be developed.

\subsection{Spatial Approximation}
To approximate the solution in space,
we use a domain decomposition method
with two non-overlapping subdomains,
where the interface is located approximately
at the center of the transition layer,
an adaptive pseudo-spectral method
on each subdomain, and collocation based
on Tchebychev polynomials, where
the collocation points are concentrated
in the transition layer.
The algorithm is standard and has been described
elsewhere~\cite{bgmm89,bbht91,p90,g94,gk97};
we summarize it here only for completeness.

Let $x^* \in (-1, 1)$
denote the (approximate) position of
the center of the transition layer;
$x^*$ varies in time ($t$),
but since $t$ enters only as a parameter
in the discussion of the spatial approximation,
we do not write it explicitly.
We decompose,
\begin{equation}
   \Omega_1 = (-1, x^*), \quad \Omega_2 = (x^*, 1) ,
   \label{partition}
\end{equation}
and map each of the subdomains
$\Omega_1$ and $\Omega_2$
linearly onto $(-1, 1)$,
\[
   g_1 :
   y \in (-1, 1)
   \mapsto
   x = g_1 (y) = - 1 + \textstyle{1 \over 2}(x^* + 1)(y + 1) \in \Omega_1 ,
\]
\[
   g_2 :
   y \in (-1, 1)
   \mapsto
   x = g_2 (y) = 1 - \textstyle{1 \over 2}(1 - x^*)(1 - y) \in \Omega_2 .
\]
The restrictions of $u$ to $\Omega_1$ and $\Omega_2$
exhibit boundary layer behavior near $x^*$.
The point $x = x^*$ corresponds to
$y = 1$ under $g_1$ and to $y = -1$ under $g_2$.
To concentrate the collocation points near $x^*$,
we define a one-parameter family of nonlinear
mappings of the interval $(-1, 1)$ onto itself,
\[
   f_1 (\cdot\,, \alpha) :
   s \in (-1, 1)
   \mapsto
   y = f_1 (s, \alpha)
   =
   1 - (4/\pi) \mbox{arctan}
   \left( \alpha \tan \textstyle{1 \over 4}(1 - s)\pi \right)
   \in (-1, 1) ,
\]
\[
   f_2 (\cdot\,, \alpha) :
   s \in (-1, 1)
   \mapsto
   y = f_2 (s, \alpha)
   =
   - 1 + (4/\pi) \mbox{arctan}
   \left( \alpha \tan \textstyle{1 \over 4}(s+1)\pi \right)
   \in (-1, 1) .
\]
If the parameter $\alpha$ is small,
$f_1$ concentrates points near $1$
and $f_2$ concentrates points near $-1$.
Concentrating points near critical points
is the computational analog of
coordinate stretching
in asymptotic analysis.
The choice of $\alpha$ can be optimized
by means of \textit{a priori}
estimates~\cite{bgmm89,bbht91};
we usually take
$\alpha = \varepsilon^{1/2}$~\cite{g94}.
The composite maps,
\[
   h_i (\cdot\,; \alpha)
   = g_i (\cdot\,) \circ f_i (\cdot\,, \alpha) ,
   \quad i=1,2 ,
\]
are one-to-one from $(-1, 1)$ onto $\Omega_i$;
we denote their inverses by
$h_i^{-1} (\cdot\,; \alpha)$, $i=1,2$.

We look for the solution of Eq.~(\ref{B})
by approximating locally
on each of the subdomains
$\Omega_1$ and $\Omega_2$
and imposing $C^1$ continuity at $x^*$.
If $U$ denotes the global approximation, then
\[
   U (x) = U_i (x)
   \quad\mbox{for } x \in \Omega_i ,
   \quad
   U \in C^1 ([-1, 1]) .
\]
The local approximations consist of
finite sums of Tchebychev polynomials,
\[
   U_i (x) =
   \sum_{j=0}^{N-1}
   a_{ij}
   T_j (h_i^{-1} (x; \alpha)) ,
   \quad
   x \in \Omega_i , \,
   i = 1, 2 ;
   \quad
   T_j (\cos\theta) = \cos (j\theta) ,
   \quad
   \theta = \pi/N .
\]

\subsection{Integration for Times of Order One}
For short-time integration it suffices
to combine a one-step forward Euler approximation
with an implicit treatment of the second-order
spatial derivative and an explicit treatment
of the nonlinear term.

Starting with an approximate solution
$U^0 = U( \cdot\,, t_0)$ at time $t_0$,
we identify the point $x^* = x^*(t_0)$ with
the location of the zero of $U^0$,
partition the domain in two subdomains,
and select the collocation points.
Fixing this configuration temporarily,
we compute a sequence
$\{U^n : n = 1, 2 \ldots\,\}$
of successive approximations $U^n$
using the algorithm
\begin{equation}
   - \varepsilon D^2 U^n
   + \frac{U^n - U^{n-1}}{\Delta t}
   + U^{n-1} DU^{n-1}
   = 0 ,
   \quad n = 1, 2, \ldots\, .
   \label{alg}
\end{equation}
The symbol $D$ represents
the pseudo-spectral differentiation operator
in physical space~\cite{chqz}.
The time step $\Delta t$ is constant,
so $U^n$ is the approximate solution of
the boundary-value problem~(\ref{B1}) at
$t_0 + n \Delta t$.
The algortithm~(\ref{alg}) is nonconservative.

When the location of the zero of $U^n$
has moved over a distance $\varepsilon$,
we suspend the algorithm~(\ref{alg}).
We shift $x^*$ to the current location of the zero,
reconfigure the partition,
update the collocation points,
replace $U^0$ by the values at the new collocation points
(using interpolation if necessary),
and continue the algorithm~(\ref{alg}).
We repeat this process until the steady state
is reached.
The change in the location of the zero
of the computed approximation becomes smaller
as time progresses,
so the same collocation configuration
serves for longer and longer time intervals.

The algorithm~(\ref{alg}) requires the solution
of $U^n$ from the equation
\begin{equation}
   A U^n = U^{n-1} + (\Delta t) U^{n-1} DU^{n-1} .
   \label{system}
\end{equation}
The matrix $A$, which is order $2N-1$,
has a block structure,
\[
   A
   =
   \pmatrix
   {A_1                       & b_1 &                            \cr
   a_1^{\mbox{\scriptsize t}} & c   & b_2^{\mbox{\scriptsize t}} \cr
                              & a_2 & A_2                        \cr} .
\]
$A_1$ and $A_2$ are square matrices of order $N-1$;
$a_1$, $a_2$, $b_1$, and $b_2$ vectors of length $N-1$;
$c$ is a constant.
The center row accounts for the
$C^1$ continuity at the interface.
This block structure allows a solution
of the system~(\ref{system})
in two parallel processes
from opposite ends.
The matrix $A$ does not change as long
as the collocation configuration is frozen.
However, its condition number increases with
the number of collocation points.
This increase puts a lower limit
on the values of $\delta$
one can handle in practice.

\begin{table}[htb]
\centering
\begin{footnotesize}
\caption{
Location of the transition layer at steady state.}
\begin{tabular}{|| c || c | c || c |  c ||} \hline\hline
         & \multicolumn{2}{c||}{$\varepsilon=0.1$}& \multicolumn{2}{c||}{$\varepsilon=0.05$} \\ \hline\hline
$\delta$ & $x^*_\infty$ & $x^*_{\rm as}$   & $x^*_\infty$ & $x^*_{\rm as}$ \\ \hline\hline
$1.0 \cdot 10^{-1}$ & 0.72464 & 0.700427 & 0.86237 & 0.850213 \\ \hline
$1.0 \cdot 10^{-2}$ & 0.47486 & 0.470176 & 0.73755 & 0.735084 \\ \hline
$1.0 \cdot 10^{-3}$ & 0.24133 & 0.240724 & 0.62055 & 0.619955 \\ \hline
$1.0 \cdot 10^{-4}$ & 0.05265 & 0.052606 & 0.50485 & 0.504826 \\ \hline
$1.0 \cdot 10^{-5}$ & 0.00537 & 0.005504 & 0.38962 & 0.389696 \\ \hline\hline
\end{tabular}
\end{footnotesize}
\end{table}
Table~1 gives the location of the transition layer
at steady state, $x^*_\infty$,
computed with the algorithm~(\ref{alg})
with $N=39$ collocation points in each subdomain
and a time step $\Delta t = 0.02$.
The number $x^*_{\rm as} = 1 - \varepsilon \ln (2/\delta)$,
which is an asymptotic estimate of $x^*$
(see Eq.~(\ref{de-as})) is given for comparison.
The initial conditions were usually obtained
by linear interpolation from the boundary data,
but variations were made to test the answers.
The lower limit on $\varepsilon$ is determined
by the fact that the computation time
for the algorithm~(\ref{alg}) increases
as $\varepsilon$ decreases.
In \S~2.4 we discuss an algorithm
suitable for long-time integration.

Table~2 shows the impact of grid refinement
(the number of collocation points, $N$)
on the value of $x^*_\infty$.
\begin{table}[htb]
\centering
\begin{footnotesize}
\caption{
Effect of grid refinement ($N$) on~$x^*_\infty$;
$\varepsilon = 0.1$,
$\delta = 1.0 \cdot 10^{-3}$.}
\begin{tabular}{|| c || c | c | c | c | c | c ||} \hline
$N$          & 15 & 19 & 29 & 39 & 49 & 59 \\ \hline
$x^*_\infty$ & 0.25576 & 0.24115 & 0.24166 & 0.24133 & 0.24140 & 0.24143 \\ \hline
\end{tabular}
\end{footnotesize}
\end{table}

\subsection{Neglecting Viscosity}
In supersensitive boundary-value problems,
the solution in the ``tails'' on either side
of the transition layer is exponentially close
to the prescribed boundary values,
so the viscous term is exponentially small there.
It is therefore tempting to assume
that one can sacrifice some accuracy
in the computation of the viscous term
during the transient phase and still find
the position of the transition layer at
steady state with a high degree of accuracy.

An extreme form of this assumption underlies
the approach where one constructs a first approximation
by ignoring the viscous term altogether.
Using a conservative finite-difference scheme,
such as Godunov, one constructs the entropy solution
of the inviscid conservation law ($\varepsilon = 0$)
and takes this as a first approximation.
One then constructs higher-order uniform approximations,
for example by means of
a heterogeneous domain-decomposition method,
using either a different numerical scheme
to solve the full viscid problem in the interior of
the transition layer
or some other approximation of
the viscous equation.

The $\chi$~method introduced by Brezzi et al.~\cite{bcr89}
is a more sophisticated nonlinear adaptive scheme
based on the same assumption.
Here, one replaces the boundary-value problem by
\begin{equation}
   - \varepsilon \chi(u_{xx}) + u_t + (f(u))_x = 0
   \mbox{ on } (-1,1) ,
   \quad
   u(-1, t) = 1 + \delta , \ 
   u(1, t) = - 1 ,
   \label{B1-chi}
\end{equation}
where $\chi \equiv \chi_{\sigma, \tau}$ is
a smooth monotone function,
$\chi(s) = 0$ if $|s| \leq \sigma$ and
$\chi(s) = s$ if $|s| \geq \sigma + \tau$
for some positive numbers $\sigma$ and $\tau$.
This method has been applied to
Burgers' equation~\cite{ac93} and
the incompressible Navier-Stokes equations~\cite{ap91}.
However, we claim that the $\chi$~method
cannot accurately predict the ultimate position
of the transition layer,
at least for Burgers' equation
on a finite interval with Dirichlet data.
This claim is supported by the following observations.

Consider the results quoted in~\cite[Table~II]{ac93}.
With few exceptions, they involve
relatively large values of $\sigma$
($\sigma$ is called $\delta$ in~\cite{ac93});
in fact, $\sigma$ is typically greater than
$\varepsilon^{-1/2}$
($\varepsilon$ is called $\nu$ in~\cite{ac93})
by one or two orders of magnitude.
The viscous term is therefore always neglected,
unless $u_{xx}$ is of the same order as $\varepsilon^{-1/2}$;
that is, the viscous term is neglected everywhere
except in the transition layer.
If the $\chi$~method gave the correct position
of the transition layer at steady state,
then the same would certainly be the case
when we simply multiply the viscous term
by a smooth function of position,
whose support is of order one and includes
the transition layer.
After all, in the latter case we account for
the viscous term over a much broader region.
These arguments lead us to consider
the boundary-value problem
\begin{equation}
   - \varepsilon H(x) u_{xx}
   + u_t 
   + u u_x
   = 0
   \mbox{ on } (-1, 1) ,
   \quad
   u (-1, t) = 1 + \delta , \;
   u (1, t) = - 1,
   \label{B1-H}
\end{equation}
instead of the boundary-value problem~(\ref{B1-chi}).
Here, $H$ is a smooth cut-off function,
\begin{equation}
   H (x)
   =
   \left\{
   \begin{array}{ll}
   \textstyle{1 \over 2} \left( 1 - \tanh (\alpha (x - x^*(t) + \beta)) \right) ,
   & \quad x < x^*(t) , \\
   \textstyle{1 \over 2} \left( 1 - \tanh (\alpha (-x + x^*(t) + \beta)) \right) ,
   & \quad x > x^*(t) .
   \end{array}
   \right.
   \label{H}
\end{equation}
We use a numerical approximation of $x^* (t)$
and choose the parameters $\alpha$ and $\beta$ so
\[
   H(x) = 1 \mbox{ if } |x-x^*(t)| < \textstyle{1 \over 2} \beta , \quad
   H(x) = 0 \mbox{ if } |x-x^*(t)| > \textstyle{3 \over 2} \beta ,
\]
to within machine accuracy ($1\cdot 10^{-15}$).
The algorithm~(\ref{alg}) leads to the solution of $U^n$
from the equation
\begin{equation}
   - \varepsilon H D^2 U^n
   +
   \frac{U^n - U^{n-1}}{\Delta t}
   + U^{n-1} D u^{n-1} 
   = 0 ,
   \quad n = 1, 2, \ldots\, .
   \label{qui}
\end{equation}
The results given in Table~3
(computed for $\varepsilon = 0.1$ and
$\delta = 1.0 \cdot 10^{-2}$,
with $\alpha = 200$)
show that the algorithm~(\ref{qui})
can give incorrect results for the position
of the transition layer at steady state.
(The correct value is $x^*_\infty = 0.47486$, see Table~1.)
\begin{table}[htb]
\centering
\begin{footnotesize}
\caption{
Location of the transition layer at steady state computed with
the $\chi$ method.}

\begin{tabular}{|| c || c | c | c | c | c | c | c ||} \hline
$\beta$      & 0.1    & 0.2    & 0.3    & 0.4    & 0.5    & 0.6    & 0.7    \\ \hline
$x^*_\infty$ & 0.0056 & 0.0110 & 0.0174 & 0.0384 & 0.0923 & 0.1067 & 0.1565 \\ \hline
\end{tabular}
\end{footnotesize}
\end{table}
The position of the shock freezes too early
during the transient phase.
The actual moment of freezing depends
on the size of the zone to the left
of the transition layer
(where $H \equiv 0$).
The result improves as $\beta$ increases,
but for $\beta = 0.8$ the algorithm fails
to converge.
(The position of the transition layer
keeps oscillating between the regions
where $H \equiv 1$ and $H \equiv 0$.)
On the other hand, if we include the missing part
of the viscous term explicitly
and use the algorithm
\begin{equation}
   - \varepsilon H D^2 U^n
   - \varepsilon (1 - H) D^2 U^{n-1}
   + \frac{U^n - U^{n-1}}{\Delta t}
   + U^{n-1} D u^{n-1}
   = 0 ,
   \quad n = 1, 2, \ldots\, ,
   \label{newqui}
\end{equation}
instead of~(\ref{qui}),
we retrieve the correct position of the transition layer
at steady state.
This result demonstrates clearly that,
when the problem is supersensitive,
it is not advisable
to neglect the viscous term,
even when that term is exponentially small.

Note that the modified algorithm~(\ref{newqui})
treats the second-order derivative explicitly
in the region where $H \equiv 0$ and
implicitly in the region where $H \equiv 1$.
The idea of using a cutoff function
to construct a composite algorithm
is described in detail in our article~\cite{gk97}.
The procedure offers a very general tool for the design
of heterogeneous domain decompositions
in the framework of a finite-difference approximation.
However, since the algorithm~(\ref{newqui})
is based on a Tchebyshev pseudo-spectral approximation,
the partially explicit treatment of the viscous term
forces a severe constraint on the time step
that cannot be circumvented.
For example,
with $\varepsilon=0.1$ and
$N = 49$ collocation points per subdomain,
the time step must be 10 times smaller
than the time step for the algorithm~(\ref{alg}).
The algorithm~(\ref{newqui})
is therefore certainly not practical
for long-time integration.

\subsection{Long-Time Integration}
The explicit treatment of the nonlinear term
in the algorithm~(\ref{alg}) imposes
a severe constraint on the time step (CFL condition).
The algorithm is therefore not suitable
for long-time integration.
The alternative approach commonly taken
is to use a fully implicit scheme
in combination with a
Newton algorithm~\cite{ao82}.
However, an implicit scheme can be expensive
and is certain to increase interprocessor
communication in a multiprocessing environment.
If the solution of the viscous conservation law
is close to a certain profile function,
as is the case for Burgers' equation,
at least after an initial transient,
the following algorithm offers a more efficient
alternative.

We start the integration of Eq.~(\ref{B})
at $t = t_0$, say, when we have an approximate
profile with a transition layer centered at $x^*$.
We construct the function $u_0$,
\begin{equation}
  u_0 (x) = - \tanh \frac{x - x^*}{2\varepsilon} ,
\end{equation}
which satifies Burgers' equation exactly,
and look for a solution $u$ of the form
\begin{equation}
  u(x, t) = u_0 (x) + \delta v(x, t) .
  \label{u-v}
\end{equation}
Then $v$ must satisfy the nonlinear boundary-value problem
\[
  - \varepsilon v_{xx}
  + v_t
  + u_0 v_x
  + u_0' v
  + \delta v v_x
  = 0
  \quad\mbox{on } (-1,1) ,
\]
\[
  v(-1, t) = \delta^{-1} (1 + \delta - u_0(-1)) , \quad
  v(1, t) = \delta^{-1} (- 1 - u_0(1)) .
\]
We integrate this boundary-value problem
for $t > t_0$, using the algorithm
\begin{equation}
  - \varepsilon D^2 V^n
  + \frac{V^n - V^{n-1}}{\Delta t}
  + u_0 D V^n
  + u_0' V^n
  =
  - \delta
  V^{n-1} D V^{n-1} ,
  \quad n = 1, 2, \ldots\,,
  \label{alg-V}
\end{equation}
and define an approximation $U$ of $u$
for $t > t_0$,
\begin{equation}
  U (x, t) = u_0 (x) + \delta V (x, t) , \quad t > t_0 .
  \label{U-u0}
\end{equation}
We proceed with the integration
as long as the supremum of
$U(\cdot\,, t) - u_0$
remains of the order of $\delta$.
When this criterion is no longer met,
at $t = t_1$ say,
we suspend the integration,
identify the point $x^*$ with
the location of the center of the
transition layer at $t_1$,
update the profile function $u_0$,
and continue the integration beyond $t_1$.
We repeat the procedure until
the steady state is reached.
Figure~1 shows a profile function $u_0$
and the computed solution $U(\cdot\,, t)$
at some time $t$.
Note that the former is monotone, the latter is not.
\begin{figure}[h!]
\begin{center}
\vspace{-5ex}
\mbox{\psfig{file=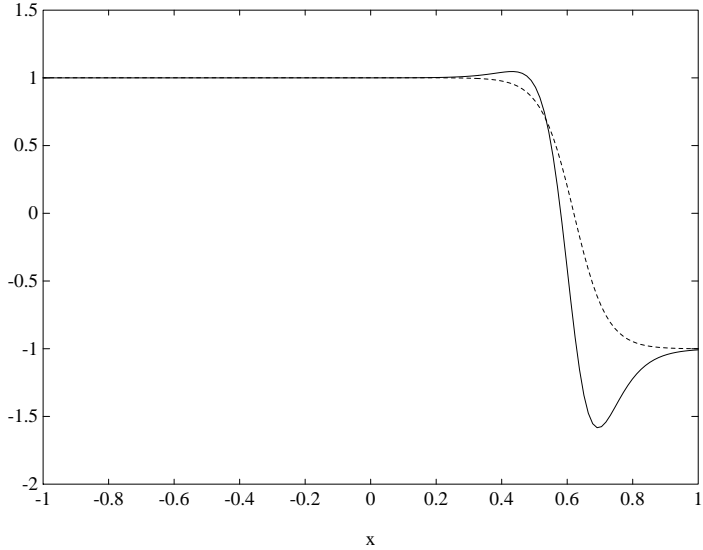,width=4.0in}}
\vspace{-2ex}
\caption{
The profile function $u_0$ (dashed line) and
the computed solution $U(\cdot\,,t)$ at some time $t$ (solid line);
$\varepsilon=0.05$, $\delta = 1 \cdot 10^{-3}$.}
\end{center}
\end{figure}

The algorithm~(\ref{alg-V}) is
very similar to~(\ref{alg}).
The spatial approximation
is handled with the same
domain-decomposition method,
so the structure of the resulting
linear system is the same as in
Eq.~(\ref{system}).
But the constraint on the time step
is relaxed by a factor $\delta$.
With the algorithm~(\ref{alg-V})
we can integrate the boundary-value problem
for values of $\varepsilon$ as small as 0.01,
down to $\delta$s of the order of $10^{-6}$,
using time steps that are typically 50 times
larger than with the algorithm~(\ref{alg}).
Table~4 gives some results for $x^*_\infty$,
computed with the algorithm~(\ref{alg-V}),
with $N=39$ collocation points in each subdomain.
\begin{table}[h!]
\centering
\begin{footnotesize}
\caption{
Location of the transition layer at steady state.}
\begin{tabular}{|| c | c | c | c | c ||} \hline\hline
                    & $\varepsilon = 0.1$ & $\varepsilon = 0.05$ & $\varepsilon = 0.02$ & $\varepsilon = 0.01$ \\ \hline\hline
$\delta$            & $x^*_\infty$ & $x^*_\infty$  & $x^*_\infty$  & $x^*_\infty$  \\ \hline\hline
$1.0 \cdot 10^{-1}$ & 0.72346      & 0.86175       &               &               \\ \hline
$1.0 \cdot 10^{-2}$ & 0.47508      & 0.73774       & 0.89518       &               \\ \hline
$1.0 \cdot 10^{-3}$ & 0.24140      & 0.62057       & 0.84827       & 0.92429       \\ \hline
$1.0 \cdot 10^{-4}$ & 0.05275      & 0.50485       & 0.80210       & 0.90104       \\ \hline
$1.0 \cdot 10^{-5}$ & 0.00561      & 0.38964       & 0.75609       & 0.87800       \\ \hline
$1.0 \cdot 10^{-6}$ & 0.00066      & 0.27452       & 0.70996       & 0.85482       \\ \hline
$1.0 \cdot 10^{-7}$ &              &               &               & 0.83084       \\ \hline\hline
\end{tabular}
\end{footnotesize}
\end{table}

\section{Two-Dimensional Case}
Next, we consider Burgers' equation
generalized to two dimensions,
\begin{equation}
   - \varepsilon \Delta u
   + u_t
   + u u_x
   + \beta u u_y
   = 0
   \quad\mbox{on } (-1, 1) \times (-\pi, \pi) ,
   \quad
   u|_{x=-1} = 1 + \delta ,
   \quad 
   u|_{x=1} = - 1 .
   \label{B-2}
\end{equation}
We assume periodicity in $y$ (period $2\pi$).
The perturbation $\delta$ may vary with $y$;
its Fourier expansion is
\begin{equation}
   \delta \equiv \delta (y)
   =
   \delta_0 \sum_{k \in {\bf Z}} \eta_k {\rm e}^{iky} .
   \label{delta-y}
\end{equation}
The coefficients $\eta_k$,
as well as the pre-factor $\delta_0$,
are independent of $y$.
The pre-factor $\delta_0$ is
arbitrarily small positive
and defined in such a way that
$\eta_0 = 1$ and
$\eta_k = O(1)$ as $\delta_0 \downarrow 0$
for $k = \pm 1, \pm 2, \ldots\,$.
The order relation~(\ref{delta-eps})
between $\delta$ and $\varepsilon$,
which must hold uniformly in $y$,
implies that
\begin{equation}
   \delta_0 = O_s ({\rm e}^{-a/\varepsilon})
   \quad\mbox{as } \delta_0, \varepsilon \downarrow 0 ,
   \label{delta0-eps}
\end{equation}
for some $a \in (0,1)$.
We show in \S~3.1 that, under these conditions,
the \textit{average} position of the transition layer
(averaged over $y$)
depends supersensitively on
the small parameters $\varepsilon$ and $\delta_0$.
In \S~3.2 we briefly review the spatial approximation.
\S~3.3 is devoted to the integration procedure.

\subsection{Profile Function}
The algorithm we propose for the solution
of the boundary-value problem~(\ref{B-2})
is again based on the assumption
that the solution is close to a known
profile function.
Our goal in this section is to show that,
under the conditions given above,
the profile function is asymptotically
independent of $y$ and again given
to leading order by the hyperbolic tangent,
as in Eq.~(\ref{tanh}).

We introduce the constant $x^* \in (0, 1)$
such that
$x^* \sim 1 - \varepsilon \ln(2/\delta_0)$
as $\varepsilon \downarrow 0$.
We define the function $u_0$,
\begin{equation}
  u_0 (x) = - \tanh \frac{x - x^*}{2\varepsilon} ,
  \label{u0}
\end{equation}
and look for a profile function $\varphi$
of the form
\begin{equation}
   \varphi (x,y) = u_0 (x) + v (x,y) .
   \label{u-0+v}
\end{equation}
Because $u_0$ satisfies Burgers' equation,
$v$ must satisfy the differential equation
\begin{equation}
  - \varepsilon \Delta v
  + u_0 v_x
  + u_0' v
  + \beta u_0 v_y
  + v v_x
  + \beta v v_y
  = 0 ,
  \label{v-eq}
\end{equation}
together with the boundary conditions
\begin{equation}
  v|_{x=-1} = \delta - \delta_0 ,
  \quad
  v|_{x=1} = 0 .
  \label{v-bc}
\end{equation}
The linearized equation,
\begin{equation}
  \ell(v)
  \equiv
  - \varepsilon \Delta v
  + u_0 v_x
  + u_0' v
  + \beta u_0 v_y
  = 0 ,
  \label{vl-eq}
\end{equation}
has a solution, $\ell(u_0') = 0$,
so we reduce the order by substituting
\begin{equation}
  v = u_0' w .
\end{equation}
Then $w$ must satisfy the equation
\begin{equation}
  - \varepsilon \Delta w
  - u_0 w_x
  + \beta u_0 w_y
  = 0 ,
  \label{w-eq}
\end{equation}
together with the boundary conditions
\begin{equation}
  w|_{x=-1} = (\delta - \delta_0)/u_0'(-1) ,
  \quad
  w|_{x=1} = 0 .
  \label{w-bc}
\end{equation}
We estimate $w$ by means of the coefficients
in its Fourier expansion,
\begin{equation}
   w (x,y) = \sum_{k \in {\bf Z}} w_k(x) {\rm e}^{iky} .
\end{equation}
The leading coefficient $w_0$ is
the solution of the boundary-value problem
\begin{equation}
  -\varepsilon w_0'' - u_0 w_0' = 0
  \quad\mbox{on } (-1,1) ,
  \quad
  w_0 (-1) = 0 ,
  \quad w_0 (1) = 0 .
  \label{w0}
\end{equation}
so $w_0 = 0$.
Note that this result is a direct consequence
of the fact that we have defined $x^*$
in terms of the average $\delta_0$
in Eq.~(\ref{u0});
any other definition leads to
an inhomogeneous boundary-value problem,
whose solution $w_0$ does not vanish.

The remaining coefficients $w_k$,
$k = \pm 1, \pm 2, \ldots\,$,
are found from the boundary-value problem
\begin{equation}
   - \varepsilon w_k''
   - u_0 w_k'
   + (\varepsilon k^2 + i\beta k u_0) w_k
   = 0 ,
   \quad w_k(-1) = \delta_0 \eta_k , \quad w_k (1) = 0 .
   \label{wk}
\end{equation}
This is a classical turning-point problem,
as $u_0$ changes sign in the interval $(-1,1)$.
The asymptotic behavior of $w_k$
as $\varepsilon \downarrow 0$
can be found by the method described
in~\cite[\S~3.E]{om91},
\begin{equation}
   w_k (x)
   \sim
   \left[
   c_k {\rm e}^{i \beta k x}
   +
   \left( 1 - c_k {\rm e}^{- i \beta k} \right) {\rm e}^{-(1+x)/\varepsilon}
   +
   \left( 0 - c_k {\rm e}^{i \beta k} \right) {\rm e}^{-(1-x)/\varepsilon}
   \right]
   w_k(-1) .
  \label{wk-as}
\end{equation}
The coefficient $c_k$ is such that the functional
\begin{equation}
   {\cal L} [w]
   =
   \frac{1}{2}
   \int_{-1}^1
   \left[
   \varepsilon w'^2
   +
   (\varepsilon k^2 + i \beta k u_0) w^2
   \right]
   \exp
   \left(
   \frac{1}{\varepsilon}
   \int_0^x u_0 (\xi) \, {\rm d}\xi
   \right) 
   \,{\rm d}x
  \label{Lw}
\end{equation}
has a critical point at $w = w_k$,
\begin{equation}
  \frac{\partial {\cal L}[w_k]}{\partial c_k} = 0 .
  \label{calL}
\end{equation}
Notice that the differential equation~(\ref{wk})
is the Euler equation of the functional ${\cal L}$;
the Neumann boundary conditions $w_k' (\pm 1) = 0$
are the natural boundary conditions associated
with ${\cal L}$.

Obviously, finding an explicit expression
for $c_k$ is out of the question.
The best we can aim for is an asymptotic expansion 
as $\varepsilon \downarrow 0$, and even here
we must resort to computational assistance.
Using the symbolic manipulation language MAPLE,
we find
\begin{equation}
  c_k \sim
  \frac
  { {\rm e}^{-(x^* + 1)/\varepsilon} }
  { (1+\beta^2) \varepsilon^2 k^2 }
  {\rm e}^{-i\beta k(1+2x^*)}
  \quad\mbox{as } \varepsilon \downarrow 0 .
  \label{ck}
\end{equation}
The first term in the brackets in Eq.~(\ref{wk-as})
represents the regular part of the asymptotic behavior of $w_k$,
which dominates in the interior;
the remaining two terms represent the singular part,
which dominates near the endpoints of the interval.
Since we are interested in the transition layer,
which is located in the interior,
we ignore the singular part
and take only the regular part,
$w_k (x) \sim c_k w_k(-1) {\rm e}^{i \beta k x}$.
That is, we take
\begin{equation}
   v_k (x)
   \sim
   \frac{\delta_0 \eta_k {\rm e}^{-(1+x^*)/\varepsilon}}
        { (1+\beta^2) \varepsilon^2 k^2} \,
   \frac{u_0'(x)}{u_0'(-1)} \,
   {\rm e}^{i\beta k(1-2x^*+x)}
   \quad\mbox{as } \varepsilon \downarrow 0 .
\end{equation}
If we use the asymptotic approximation
\[
   \frac{u_0'(x)}{u_0'(-1)}
   =
   \frac
      {1 - \tanh^2 ((x-x^*)/(2\varepsilon))}
      {1 - \tanh^2 ((-1-x^*)/(2\varepsilon))}
   \sim
   \frac
   {1 - \tanh^2 ((x-x^*)/(2\varepsilon))}
   {4{\rm e}^{-(1+x^*)/\varepsilon}} ,
\]
we obtain the asymptotic expression
\begin{equation}
   v_k (x)
   \sim
   \frac{ \delta_0 \eta_k }
        { 4 (1+\beta^2) \varepsilon^2 k^2} \,
   {\rm e}^{i\beta k(1-2x^*+x)}
   \left(
   1 - \tanh^2 \frac{x-x^*}{2\varepsilon}
   \right) .
  \label{vk-as}
\end{equation}
This result implies that the Fourier series
of $v$, as well as those of $v_x$ and $v_y$,
converge.
Furthermore, $\|v\|_\infty$ and
$\|v_y\|_\infty$ are
$O(\delta_0 \varepsilon^{-2})$.

Finding the asymptotic behavior of
$\|v_x\|_\infty$ is less obvious.
It follows from Eq.~(\ref{wk-as}) that
$w'_k (x) \sim c_k {\rm e}^{i\beta kx}$
as $\varepsilon \downarrow 0$, at least for
$x$ in the interior of the domain.
Therefore,
$v'_k(x) = (u_0' w_k)' (x)
\sim (u_0'' w_k)(x)
\sim \varepsilon^{-1} (u_0' w_k)(x)
= \varepsilon^{-1} v_k(x)$.
Hence, 
$\|v_x\|_\infty = O(\delta_0 \varepsilon^{-3})$
as $\delta_0, \varepsilon \downarrow 0$.

Since $\delta_0 = O_s ({\rm e}^{-a/\varepsilon} )$
for some $a \in (0,1)$,
we have
$\delta_0 \varepsilon^{-p}
= O_s ({\rm e}^{-a'/\varepsilon})$
($p = 2, 3$) for any $a' \in (0,a)$,
so any solution $v$ of Eq.~(\ref{vl-eq})
which satisfies the boundary conditions~(\ref{v-bc})
is transcendentally small.
The residue $v v_x + \beta v v_y$,
which was ignored in the transition
from the nonlinear equation~(\ref{v-eq})
to the linear equation~(\ref{vl-eq})
is likewise transcendentally small and,
in fact, $O(\delta_0^2 \varepsilon^{-5})$,
so we also have an a posteriori
justification for the linearization.

These arguments motivate the choice
of $u_0$, which depends only on $x$,
as the profile function in the design
of the numerical algorithms for Eq.~(\ref{B-2}).

\subsection{Spatial Approximation}
The spatial approximation is again based
on a domain decomposition with two
non-overlapping subdomains on either side
of the $y$-averaged location
of the center of the transition layer.
On each subdomain we use an adaptive
pseudo-spectral method in the $x$ direction
and a finite-difference method in the $y$ direction.
The pseudo-spectral method is the same as in the
one-dimensional case;
it uses Tchebychev polynomial collocation
with $N_x$ collocation points.

Since the transition layer is close
to a plane parallel to the $x$ axis, 
there is no need to resort to
an adaptive grid in the $y$ direction.
For our numerical experiments
we chose a regular grid with mesh width
$h = 2 \pi / N_y$ and
a sixth-order central finite-difference
approximation of $u_{yy}$ and $u_y$.
The choice may seem inconsistent
with the spectral approximation
in the $x$ direction;
a more obvious choice would be
a Fourier approximation
in the $y$ direction.
Theoretically, the finite-difference approximation
in the $y$ direction restricts the accuracy
of the approximation for a regular problem
($\varepsilon = O_s(1)$)
to sixth order,
less than the accuracy guaranteed by
the pseudo-spectral approximation
in the $x$ direction.
However, as the transition layer is close
to a plane parallel to the $x$ axis,
it is relatively easy to keep the numerical error
in the finite-difference approximation of the term
$\varepsilon u_{yy}$ smaller than
the numerical error in the pseudo-spectral approximation
of the term $\varepsilon u_{xx}$
with a moderate number of discretization
points $N_y$.
There is, therefore, no need to use
a better approximation,
like the Fourier approximation,
for the term $\varepsilon u_{yy}$.
Furthermore,
the spectral radius of $D_y^2$
is smaller with sixth-order finite differences
than with Fourier differentiation.
This difference implies an additional
advantage for a finite-difference approximation
when the $y$ derivatives are treated explicitly~\cite{dg95}.

\subsection{Integration for Times of Order One}
We extend the Euler scheme~(\ref{alg})
to two dimensions as follows:
\begin{equation}
  - \varepsilon D^2_x U^n
  + \frac{U^n - U^{n-1}}{\Delta t}
  =
  \varepsilon D^2_y U^{n-1}
  + U^{n-1} D_x U^{n-1}
  + \beta U^{n-1} D_y U^{n-1} ,
  \quad n = 1, 2, \ldots \,.
  \label{2dalg}
\end{equation}
Here, $D_x$ is the pseudo-spectral differential operator
with Tchebychev polynomials,
$D_y$ the finite-difference operator
with sixth-order central finite differences.

The algorithm~(\ref{2dalg}) has
several features that make it
readily parallelizable.
First, the approximations
$D_y U^{n-1}$ and $D_y^2 U^{n-1}$
are taken explicitly, so
the variable $y$ is only a parameter.
Second, because we are using
finite-difference approximations,
we have only local data dependencies.
This latter point especially offers
a significant advantage over a spectral
method, which uses global interpolation.

Table~5 shows the results for
the boundary-value problem~(\ref{B-2})
with $\beta = 1$,
$\varepsilon = 0.1$
and piecewise constant boundary data
with $\delta_0 = 1.0 \cdot 10^{-2}$,
\begin{equation}
  u(-1,y) = \left\{
  \begin{array}{ll}
  1.01 - \Delta\delta & \mbox{if } -\pi \leq y < -\textstyle{1 \over 2}\pi ,    \\
  1.01 + \Delta\delta & \mbox{if } -\textstyle{1 \over 2}\pi \leq y < \textstyle{1 \over 2}\pi ,\\
  1.01 - \Delta\delta & \mbox{if } \textstyle{1 \over 2}\pi \leq y < \pi .
  \end{array}
  \right.
  \label{pcbc}
\end{equation}
The algorithm~(\ref{2dalg}) was applied with
$N_x = 39$ collocation points per subdomain
in the $x$ direction
and $N_y = 32$ interpolation points in the $y$ direction.
The table gives the $y$-averaged location of the
transition layer at steady state, $<\!x^*_\infty\!>$,
as well as the maximum deviation
of the center of the transition layer
from its $y$-averaged location, $\Delta x^*$;
that is, $x^*_\infty (y)$ varies between
$<\!x^*_\infty\!> - \, \Delta x^*$ and
$<\!x^*_\infty\!> + \, \Delta x^*$.
\begin{table}[htb]
\centering
\begin{footnotesize}
\caption{
Location of the transition layer at steady state; boundary data~(\ref{pcbc}).}
\begin{tabular}{|| c | c | c ||} \hline\hline
$\Delta\delta$        & $<\!x^*_\infty\!>$ & $\Delta x^*$           \\ \hline\hline
$0.25 \cdot 10^{-2}$  & 0.4758         & $1.3114 \cdot 10^{-2}$ \\ \hline
$0.50 \cdot 10^{-2}$  & 0.4759         & $2.3584 \cdot 10^{-2}$ \\ \hline
$1.0 \cdot 10^{-2}$   & 0.4758         & $4.8865 \cdot 10^{-2}$ \\ \hline
$1.5 \cdot 10^{-2}$   & 0.4750         & $7.8632 \cdot 10^{-2}$ \\ \hline
$2.0 \cdot 10^{-2}$   & 0.4738         & $9.3250 \cdot 10^{-2}$ \\ \hline
$3.0 \cdot 10^{-2}$   & 0.4700         & $15.688 \cdot 10^{-2}$ \\ \hline\hline
\end{tabular}
\end{footnotesize}
\end{table}
These results show
that the algorithm~(\ref{2dalg})
is extremely effective for the boundary-value problem.
As $\Delta\delta$ increases,
$\Delta x^*$ grows approximately linearly
with~$\Delta\delta$.
The graph of $U - u_0$ maintains its overall shape,
so its width (which measures $\Delta x^*$)
varies in proportion to its height
(which measures $\| U - u_0 \|_\infty$).
The numerical results therefore indicate
also that $\| U - u_0 \|_\infty$
grows approximately linearly
with~$\Delta\delta$.
This conclusion matches the results
of the asymptotic analysis in \S~3.1,
in particular Eq.~(\ref{vk-as}),
where it was shown that
$v_k$ is proportional to $\eta_k$.

Results for a much harder case are presented
in Table~6.
The parameters $\beta$ and $\varepsilon$ are fixed
as before,
$\beta = 1$ and $\varepsilon = 0.1$,
but this time the data at the left boundary
are sharply peaked at the midpoint,
\begin{equation}
  u(-1,y)
  =
  1.005 + (\Delta\delta) {\rm e}^{-20(1 - \cos y)} ,
  \quad -\pi < y < \pi .
  \label{bpbc}
\end{equation}
The algorithm~(\ref{2dalg}) was again applied with
$N_x = 39$ collocation points per subdomain
in the $x$ direction
and $N_y = 32$ interpolation points in the $y$ direction.
The table gives, in addition to the values of
$<\!x^*_\infty\!>$ and $\Delta x^*$,
the asymptotic value
$x^*_{\rm as} = 1 - \varepsilon \ln (2/\delta_0)$.
\begin{table}[htb]
\centering
\begin{footnotesize}
\caption{
Location of the transition layer at steady state; boundary data~(\ref{bpbc}).}
\begin{tabular}{|| c | c | c | c ||} \hline\hline
$\Delta\delta$      &$<\!x^*_\infty\!>$&$x^*_{\rm as}$&$\Delta x^*$ \\ \hline\hline
$0.25 \cdot 10^{-2}$& 0.40808 & 0.40811 & $0.27278 \cdot 10^{-2}$ \\ \hline
$0.50 \cdot 10^{-2}$& 0.41277 & 0.41241 & $0.37231 \cdot 10^{-2}$ \\ \hline
$1.0 \cdot 10^{-2}$ & 0.42096 & 0.42052 & $0.92263 \cdot 10^{-2}$ \\ \hline
$2.0 \cdot 10^{-2}$ & 0.43573 & 0.43506 & $1.4589 \cdot 10^{-2}$ \\ \hline
$3.0 \cdot 10^{-2}$ & 0.44854 & 0.44782 & $2.3924 \cdot 10^{-2}$ \\ \hline\hline
\end{tabular}
\end{footnotesize}
\end{table}
The average location of the transition layer is predicted
very well by the asymptotics.
Again, the maximum deviation $\Delta x^*$
grows with $\Delta\delta$,
although not linearly as in the case of
the step boundary data~(\ref{pcbc}).

Table~7 shows the impact of grid refinement
(the number of collocation points per subdomain, $N_x$,
and the number of discretization points, $N_y$)
on the computed value of $x^*_\infty (y)$
for the problem with boundary data~(\ref{pcbc}),
$\Delta\delta = 0.01$.
\begin{table}[htb]
\centering
\begin{footnotesize}
\caption{
Effect of grid refinement ($N_x$, $N_y$)
on $<\!x^*_\infty\!>$ (upper entries)
and $\Delta x^*$ (lower entries);
boundary data~(\ref{pcbc}) with $\Delta\delta = 0.01$.}
\begin{tabular}{|| c || c | c | c | c | c ||} \hline\hline
$N_y$ & $N_x=19$ & $N_x=29$ & $N_x=39$ & $N_x=49$ & $N_x=59$         \\ \hline\hline
8     & 0.47848  & 0.47653  & 0.47523  & 0.47557  & 0.47547 \\
      & 4.7396 $10^{-2}$    & 4.8375 $10^{-2}$    & 4.8902 $10^{-2}$ 
      & 5.2615 $10^{-2}$    & 5.2251 $10^{-2}$              \\ \hline
16    & 0.47807  & 0.47614  & 0.47610  & 0.47478  & 0.47516 \\
      & 4.7396 $10^{-2}$    & 5.8062 $10^{-2}$    & 4.8849 $10^{-2}$
      & 4.9273 $10^{-2}$    & 4.9518 $10^{-2}$              \\ \hline
32    & 0.47883  & 0.47611  & 0.47578  & 0.47536  & 0.47481 \\
      & 4.7396 $10^{-2}$    & 4.8374 $10^{-2}$    & 4.8864 $10^{-2}$
      & 4.9267 $10^{-2}$    & 4.9516 $10^{-2}$              \\ \hline
64    & 0.47847  & 0.47610  & 0.47561  & 0.47513  & 0.47524 \\
      & 4.7396 $10^{-2}$    & 4.8372 $10^{-2}$    & 4.8867 $10^{-2}$
      & 4.9266 $10^{-2}$    & 5.2280 $10^{-2}$              \\ \hline\hline
\end{tabular}
\end{footnotesize}
\end{table}

An obvious way to parallelize the algorithm~(\ref{2dalg})
is to partition the interval $(-\pi, \pi)$
into subintervals of equal length $2\pi/N_y$
and map this partition onto a ring of processors.
Thus, one can achieve high speedups
on a Paragon using nonblocking communications.
If each processor covers at least
four mesh points in the $y$ direction,
only nearest-neighbor communication
is needed.
Table~8 gives ample evidence
that the algorithm~(\ref{2dalg})
is highly scalable;
doubling the number of processors
with the problem size results in
almost identical CPU times.
\begin{table}[htb]
\centering
\begin{footnotesize}
\caption{CPU time for 1,000 time steps on the Paragon XP/S
as a function of the number of processors ($P$)
and the size of the problem (measured by $N_y$);
$N_x = 49$.}
\begin{tabular}{|| c || c | c | c | c | c ||} \hline\hline
$N_y$&  $P=1$ &  $P=2$ &  $P=4$ & $P=8$  & $P=16$\\ \hline\hline
32   & 129.16 &  65.88 &  33.91 &  17.43 &       \\ \hline
64   & 252.85 & 130.02 &  65.85 &  33.85 & 17.45 \\ \hline
128  & 519.05 & 257.85 & 129.68 &  65.85 & 33.92 \\ \hline
256  &        & 515.01 & 257.49 & 129.77 & 65.93 \\ \hline\hline
\end{tabular}
\end{footnotesize}
\end{table}

Additional parallelism can be introduced
by decomposing the domain in the $x$ direction.
However, our experience with a similar algorithm
for combustion problems indicates a potentially
significant decrease (as much as 70\%)
in the efficiency of the algorithm~\cite{gtd97}.

In general, the algorithm~(\ref{2dalg})
is very well adapted to
the quasi one-dimensional structure
of the transition layer.
The algorithm predicts
the location of the transition layer
at steady state with a significant accuracy.
The time step is of the same order of magnitude
as for the one-dimensional analog~(\ref{alg}).

\subsection{Long-Time Integration}
The algorithm~(\ref{2dalg})
needs to be modified
for long-time integration.
We distinguish between the cases
$\beta = 0$ and $\beta \not= 0$.

If $\beta = 0$, we use an algorithm
similar to the one described in \S~2.4.
We start the integration of Eq.~(\ref{B-2})
at $t = t_0$, say.
We identify the point $x^*$ with the location
of the zero of the approximation $U$ of $u$,
averaged over $y$, at $t = t_0$
and define the profile function
$u_0$ as in Eq.~(\ref{u0}),
\begin{equation}
  u_0 (x)
  = - \tanh \frac{x-x^*}{2\varepsilon} ,
  \label{2du0}
\end{equation}
Then we integrate the nonlinear boundary-value problem
\[
  - \varepsilon \Delta v
  + v_t
  + u_0 v_x
  + u_0' v
  + \delta_0 \varepsilon^{-2} v v_x
  = 0
  \quad\mbox{on } (-1,1) \times (-\pi,\pi)
\]
forward in time,
subject to the boundary conditions
\[
  v|_{x=-1} = \delta_0^{-1} \varepsilon^2 (1 + \delta - u_0(-1)) , \quad
  v|_{x=1} = \delta_0^{-1} \varepsilon^2 (- 1 - u_0(1)) ,
\]
using the algorithm
\begin{equation}
  - \varepsilon D_x^2 V^n
  + \frac{V^n - V^{n-1}}{\Delta t}
  + u_0 D_x V^n
  + u_0' V^n
  =
  \varepsilon D_y^2 V^{n-1}
  - \delta_0 \varepsilon^{-2}
  V^{n-1} D_x V^{n-1} ,
  \label{2dalg-V}
\end{equation}
for $n = 1, 2, \ldots\,$.
We define the approximation $U$ of $u$,
\begin{equation}
  U (x,y,t) = u_0 (x) + \delta_0 \varepsilon^{-2} V (x,y,t) ,
  \label{2dU-u0}
\end{equation}
and integrate as long as the supremum of
$U(\cdot\,,\cdot\,, t) - u_0$
remains of the order of $\delta_0 \varepsilon^{-2}$.
When this criterion is no longer met,
at $t = t_1$ say,
we suspend the integration,
identify the point $x^*$ with
the location of the center of the transition layer
(averaged over $y$),
and update the profile function $u_0$.
We repeat the procedure until
the steady state is reached.

The time step for the algorithm~(\ref{2dalg-V})
is limited by the (explicit) term
$\varepsilon D_y^2 V^{n-1}$,
$\Delta t < c (2\pi/N_y)^2 / \varepsilon$,
for some constant $c < \textstyle{1 \over 2}$.
This limitation is not too severe,
as $\varepsilon$ is very small and 
the variation of the solution in
the $y$ direction is exponentially small,
so $N_y$ need not be large.

If $\beta \not= 0$, the situation becomes
more complicated.
One can, of course, extend
the algorithm~(\ref{2dalg-V})
trivially by incorporating
the convective term
in the right member,
\[ 
  - \varepsilon D_x^2 V^n
  + \frac{V^n - V^{n-1}}{\Delta t}
  + u_0 D_x V^n
  + u_0' V^n
\]
\begin{equation}
  =
  \varepsilon D_y^2 V^{n-1}
  - \delta_0 \varepsilon^{-2}
  V^{n-1} D_x V^{n-1}
  - \beta u_0 D_y V^{n-1}
  - \beta \delta_0 \varepsilon^{-2}
  V^{n-1} D_y V^{n-1} ,
  \label{2dalg-V-beta}
\end{equation}
for $n = 1, 2, \ldots\,$.
Representative results
obtained in this way
for the boundary-value problem~(\ref{B-2})
with $\beta = 1$ and $\varepsilon = 0.02$
are given in Table~9.
The boundary data are again piecewise constant,
as in Table~5, but with $\delta_0 = 1.0 \cdot 10^{-6}$,
\begin{equation}
  u(-1,y) = \left\{
  \begin{array}{ll}
  1 + 1.0 \cdot 10^{-6} - \Delta\delta & \mbox{if } -\pi \leq y < -\textstyle{1 \over 2}\pi , \\
  1 + 1.0 \cdot 10^{-6} + \Delta\delta & \mbox{if } -\textstyle{1 \over 2}\pi \leq y < \textstyle{1 \over 2}\pi ,\\
  1 + 1.0 \cdot 10^{-6} - \Delta\delta & \mbox{if } \textstyle{1 \over 2}\pi \leq y < \pi .
  \end{array}
  \right.
  \label{pcbc-beta}
\end{equation}
The algorithm~(\ref{2dalg-V-beta}) was applied with
$N_x = 39$ collocation points per subdomain
in the $x$ direction
and only $N_y = 16$ interpolation points
in the $y$ direction.
\begin{table}[htb]
\centering
\begin{footnotesize}
\caption{
Location of the transition layer at steady state;
boundary data~(\ref{pcbc-beta}).}
\begin{tabular}{|| c | c | c ||} \hline\hline
$\Delta\delta$       & $<\!x^*_\infty\!>$ & $\Delta x^*$           \\ \hline\hline
$1.0 \cdot 10^{-6}$  & 0.7099         & $<1.0 \cdot 10^{-4}$ \\ \hline
$1.0 \cdot 10^{-5}$  & 0.7099         & $<1.0 \cdot 10^{-4}$ \\ \hline
$1.0 \cdot 10^{-4}$  & 0.7101         & $3.9 \cdot 10^{-3}$ \\ \hline
$0.5 \cdot 10^{-3}$  & 0.7122         & $1.9 \cdot 10^{-2}$ \\ \hline
$1.0 \cdot 10^{-3}$  & 0.7162         & $3.0 \cdot 10^{-2}$ \\ \hline\hline
\end{tabular}
\end{footnotesize}
\end{table}
We observe that the average location
of the transition layer does not change
appreciably as long as $\Delta\delta$
is of the same order as the average perturbation $\delta_0$.
As $\Delta\delta$ increases,
the perturbation is no longer small
compared with $\delta_0$,
and the asymptotic results of \S~3.1
do not necessarily apply.
Indeed, $\Delta x^*$ does not appear
to vary linearly with $\Delta\delta$,
as was the case in Table~6.

\begin{figure}[h!]
\begin{center}
\mbox{\psfig{file=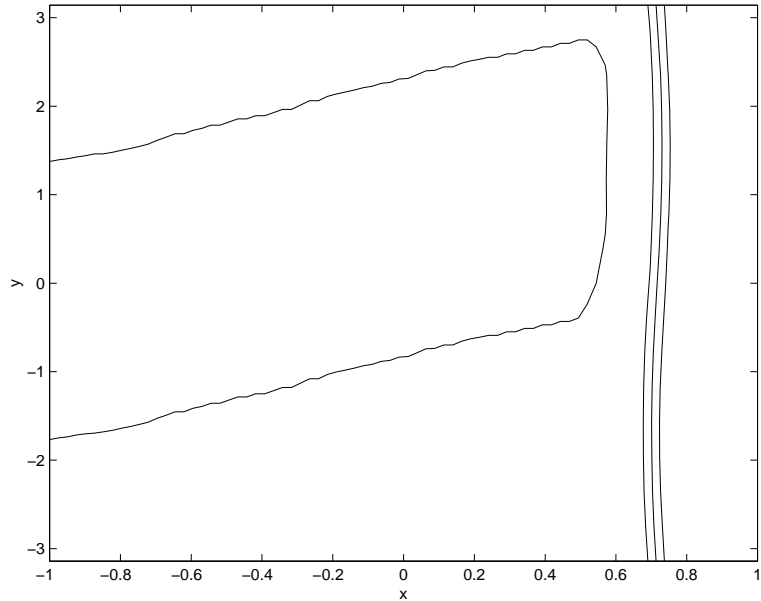,width=4.0in}}
\vspace{-2ex}
\caption{
Contour lines of the solution $U$ at steady state.}
\end{center}
\end{figure}
\begin{figure}[h!]
\begin{center}
\mbox{\psfig{file=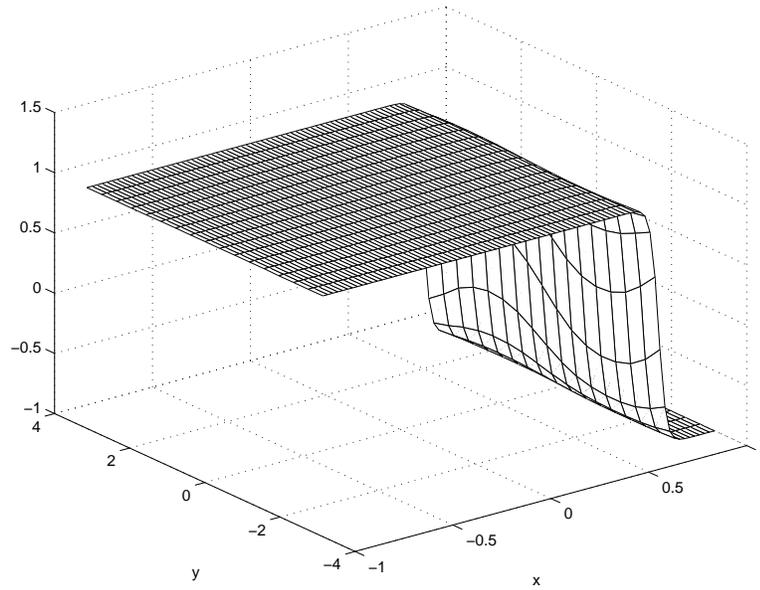,width=4.0in}}
\vspace{-2ex}
\caption{
The solution $U$ at steady state.}
\end{center}
\end{figure}
\begin{figure}[h!]
\begin{center}
\mbox{\psfig{file=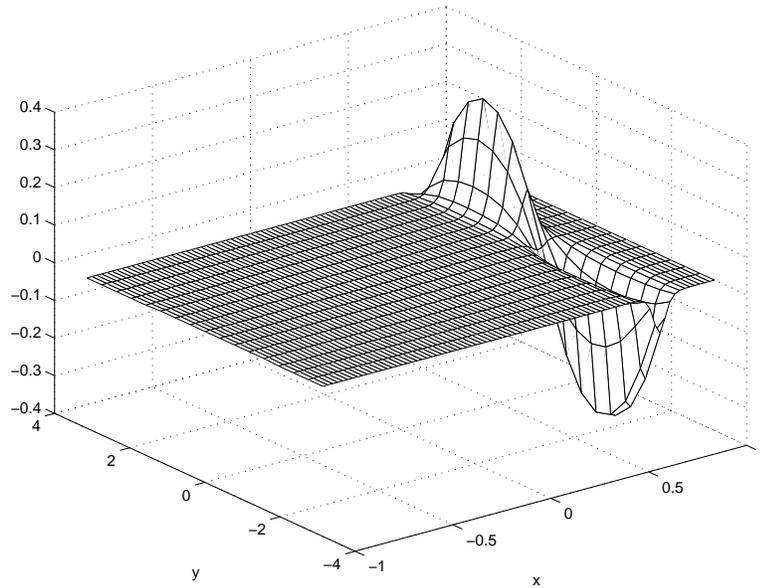,width=4.0in}}
\vspace{-2ex}
\caption{
The difference $U - u_0$ at steady state.}
\end{center}
\end{figure}

A graphical representation of the computed solution $U$
and its deviation from the profile function,
$U - u_0$, for the case
$\Delta \delta = 1.0 \cdot 10^{-3}$
are given in Fig.~2
(contour lines of $U$ at steady state)
and Figs.~3 and~4
(perspective drawings of $U$ and
the difference $U - u_0$ at steady state).
The data for these figures were taken after
2,500,000 time steps ($\Delta t = 0.4$).

The results are quite good,
but for long-time integration
one would do better
by looking for $v$ in terms of
its Fourier coefficients $v_k$
and adopting an implicit scheme
in the $y$ direction.
Thus, the constraint on the time step
becomes the same as in the one-dimensional case.
The algorithm would be of the following type:
\begin{equation}
  - \varepsilon D_x^2 V_k^n
  + \frac{V_k^n - V_k^{n-1}}{\Delta t}
  + u_0 D_x V_k^n
  + (i\beta k u_0 + \varepsilon k^2) V_k^n
  + u_0' V^n
  =
  - \delta_0 \varepsilon^{-2}
  W_k^{n-1} ,
  \label{2dalg-Vbeta}
\end{equation}
for $n = 1, 2 , \ldots\,$;
$D_y$ is the matrix of differentiation
with respect to $y$ in Fourier space,
and $W_k$ is some approximation
to the $k$th Fourier coefficient of
$vv_x + \beta vv_y$.

This algorithm parallelizes with respect
to the Fourier modes.
It has been applied to a problem involving
a propagating combustion front
in a moving fluid~\cite{gtd97,gtd98}.

\end{document}